\documentclass{amsart}
\usepackage{amsmath,amssymb,amscd,latexsym}
\usepackage[all]{xy}

\newcommand{\dist}{\mathrm{dist}}

\newcommand{\supp}{\mathrm{supp }\,}

\newcounter{number}[section]

\newenvironment{nummer}{\refstepcounter{number}{\noindent\arabic{section}.\arabic{number}}}{}

\newcommand{\bn}{\noindent \begin{nummer} \rm}
\newcommand{\en}{\end{nummer}}

\newenvironment{ntheorem}{\noindent {\sc Theorem:} \it}{}
\newenvironment{nlemma}{\noindent {\sc Lemma:} \it}{}
\newenvironment{nprop}{\noindent {\sc Proposition:} \it}{}

\newenvironment{ncor}{\noindent {\sc Corollary:} \it}{}

\newenvironment{nproof}{\noindent {\sc Proof:}}{\mbox{}\hfill 
\rule[-.2ex]{.25em}{1.8ex}}

\begin{document}

\title[Classifying crossed products]{{\sc Classifying crossed product $\mathrm{C}^{*}$-algebras}}

\author{Wilhelm Winter}
\address{Mathematisches Institut\\
Universit\"at M\"unster\\ Germany}

\email{wwinter@uni-muenster.de}

\date{\today}
\thanks{Supported by EPSRC  Grant EP/I019227/1, DFG SFB 878, GIF Grant 1137-30.6/2011}

\setcounter{section}{-1}

\begin{abstract}
I combine recent results in the structure theory of nuclear $\mathrm{C}^{*}$-algebras and in topological dynamics to classify certain types of crossed products in terms of their Elliott invariants. 
In particular, transformation group $\mathrm{C}^{*}$-algebras associated to free minimal $\mathbb{Z}^{d}$-actions on the Cantor set with compact space of ergodic measures are classified by their ordered $\mathrm{K}$-theory. In fact, the respective statement holds for finite dimensional compact metrizable spaces, provided that projections of the crossed products separate tracial states.
Moreover, $\mathrm{C}^{*}$-algebras associated to certain minimal homeomorphisms of spheres $S^{2n+1}$ are only determined by their spaces of invariant Borel probability measures (without a condition on the space of ergodic measures). 
Finally, I show that for a large collection of classifiable $\mathrm{C}^{*}$-algebras, crossed products by $\mathbb{Z}^{d}$-actions are generically again classifiable.
\end{abstract}

\maketitle

\section{Introduction}

\noindent
The structure and classification theory of nuclear $\mathrm{C}^{*}$-algebras has seen substantial progress in recent years, largely spurred by the interplay of certain topological, algebraic and homological regularity properties. These allow for some amount of interpretation, but generally arise as finite topological dimension, tensorial absorption of certain touchstone $\mathrm{C}^{*}$-algebras and order completeness of homological invariants, cf.\ \cite{EllToms:BullAMS} and \cite{Winter:dr-Z-stable}. 

When Andrew Toms and I began to study the tight connections between finite decomposition rank (cf.\ \cite{KirWinter:dr}), $\mathcal{Z}$-stability (cf.\ \cite{JiaSu:Z,RorWin:Z-revisited}) and strict comparison (cf.\ \cite{Ror:Z-absorbing}), this was for very specific (and, to some extent artificial) classes of simple inductive limits (see \cite{TomsWinter:VI}). However, the three properties soon turned out to occur (or to fail) simultaneously in much broader generality. We then conjectured that they are equivalent for all separable, simple, unital, nonelementary, nuclear $\mathrm{C}^{*}$-algebras which are finite (hence admit a tracial state). With the introduction of nuclear dimension in \cite{WinterZac:dimnuc}, this point of view -- and the regularity conjecture -- became available also in a not necessarily finite setting. 

Of the implications in the regularity conjecture, the problem of when $\mathcal{Z}$-stability implies finite nuclear dimension or even finite decomposition rank arguably remains the most intriguing. While by now there are several results establishing such implications, most of these factorize through classification theorems of some sort. In \cite{TikWin:Z-dr}, a more direct argument (based on \cite{KirRor:pi3}) was given; this does not use a classification result in any way and in fact it even works in a not necessarily simple situation, but it requires the algebra to be locally homogeneous -- a rather strong structural hypothesis.  

In \cite{MatSat:dr-UHF}, Matui and Sato obtained a very convincing result in the simple, nuclear, monotracial case: In this situation, $\mathcal{Z}$-stability implies finite decomposition rank provided the algebra is also quasidiagonal. They moreover showed that in the traceless case, $\mathcal{Z}$-stability implies finite nuclear dimension. These results are in line with Kirchberg--Phillips classification of purely infinite $\mathrm{C}^{*}$-algebras (the traceless, finite nuclear dimension case), and with Lin's classification of TAF algebras (these have traces, and are always quasidiagonal). Of course, it remains an important open question whether tracial $\mathrm{C^{*}}$-algebras with finite nuclear dimension are automatically quasidiagonal. 

In this paper, I am interested in applying these ideas to the classification problem for simple, nuclear $\mathrm{C}^{*}$-algebras, following Elliott's program to classify nuclear $\mathrm{C}^{*}$-algebras by $\mathrm{K}$-theoretic invariants. The problem is particularly relevant for so-called transformation group $\mathrm{C}^{*}$-algebras, i.e., $\mathrm{C}^{*}$-algebras associated to topological dynamical systems via the crossed product construction. These are the perhaps most natural source of examples, and they are important invariants of dynamical systems as has been impressively demonstrated in \cite{GPS:orbit} (for Cantor minimal $\mathbb{Z}$-actions). Also, for such crossed products there are highly useful tools for computing the invariant, in particular the Pimsner--Voiculescu exact sequence.

For minimal actions on finite dimensional compact spaces, their $\mathrm{C}^{*}$-algebras were successfully classified by ordered $\mathrm{K}$-theory  in \cite{TomsWinter:minhom} (see also \cite{TomsWinter:PNAS}), at least if projections separate tracial states (e.g., in the uniquely ergodic situation). While this is already quite satisfactory, it begs to be generalized in two directions: 

First, what if projections do not separate tracial states? This is interesting even when there are no nontrivial projections and only finitely many extremal traces -- and such examples indeed do exist, like $\mathrm{C}^{*}$-algebras of minimal homeomorphisms of odd dimensional spheres, as considered in \cite{Con:Thom} and \cite{Windsor:finite-ergodic} (those examples rely on the fast approximation method of \cite{FatHer:diffmin}, first introduced by Anosov and Katok). 

Second, what about more general group actions? Even $\mathbb{Z}^{d}$-actions on Cantor sets are notoriously difficult to handle in this respect (but cf.\ \cite{GMPS:orbit-equivalence-Z2, GMPS:orbit-equivalence-Zd}). One of the reasons is that, on the $\mathrm{C}^{*}$-side, analogues of Putnam's orbit breaking subalgebras are not easy to find, let alone to use.  

Below I combine recent results on the structure of dynamical systems and crossed products with a new technique (see Theorem~\ref{A}) to make progress on both of these questions. This new method reduces the problem of showing that a $\mathrm{C}^{*}$-algebra is classifiable to the problem of embedding it into a classifiable $\mathrm{C}^{*}$-algebra in a sufficiently nice fashion. (Here, by `classifiable' I mean TAF, or TAI, in the sense of Lin -- for the crossed products in question we do not have to worry about the Universal Coefficient Theorem.) I am confident that this method will prove to be as useful for classification as the one introduced in \cite{Win:localizingEC} (where classification up to $\mathcal{Z}$-stability was reduced to classification up to UHF-stability), see also \cite{LinNiu:KKlifting, Lin:asu-class}.

Apart from the aforementioned Theorem~\ref{A}, the main ingredients are a recent result of Szab\'o, where he establishes finite Rokhlin dimension (in the sense of \cite{HirWinZac:Rokhlin-dimension}) of free $\mathbb{Z}^{d}$-actions of finite dimensional spaces, a recent result of Strung which ensures that certain transformation group $\mathrm{C}^{*}$-algebras nicely embed into classifiable models, and a result of Lin which establishes a strong form of AF-embeddability of certain types of crossed product $\mathrm{C}^{*}$-algebras. 

As byproducts I obtain a simple proof of a special case of \cite[Theorem~6.1]{MatSat:dr-UHF} and show that, for certain classifiable $\mathrm{C}^{*}$-algebras, generic sets of $\mathbb{Z}^{d}$-actions give rise to crossed products which are again classifiable.

Section~{\ref{section1}} recalls a characterization of TA$\mathcal{S}$ algebras (with $\mathcal{S}$ a suitable class of building blocks) which will be useful for our purposes. It also summarizes some of the relevant results on the classification of rationally TAF (or TAI) $\mathrm{C}^{*}$-algebras. Section~{\ref{section2}} contains the main technical result, along with a Corollary which illustrates the `classification by embeddings' method in a somewhat more concise (and less technical) manner. We also obtain a short proof of a special case of a recent result of Matui and Sato along these lines. In Section~{\ref{section3}} Theorem~\ref{A} is applied to transformation group $\mathrm{C}^{*}$-algebras of free and minimal $\mathbb{Z}^{d}$-actions.  The result in particular covers free and minimal Cantor $\mathbb{Z}^{d}$-actions with compact space of ergodic measures; the $\mathrm{C}^{*}$-algebras in this case are classified by their ordered $\mathrm{K}$-theory. In Section~{\ref{section4}} it is shown that $\mathrm{C}^{*}$-algebras associated to minimal homeomorphisms of spheres (of odd dimension at least 3) obtained by the fast approximation method of \cite{FatHer:diffmin} are classified by their spaces of invariant Borel probability measures. Finally, Section~{\ref{section5}} shows that (in a suitable context) classifiability generically passes to crossed products by $\mathbb{Z}^{d}$-actions.

I would like to thank the referee for carefully proofreading the paper and for a number of helpful comments and suggestions.

\section{TA$\mathcal{S}$ algebras and classification}
\label{section1}

\bn
\label{M}
For convenience, let us recall the following characterization of  rationally TA$\mathcal{S}$ algebras from \cite[Lemma~1.2]{StrWin:UHFslicing} (see also \cite[Lemma~3.2]{Win:Z-class}). By $\mathcal{Q}$ we denote the universal UHF $\mathrm{C}^{*}$-algebra.

\begin{nprop}
Let $\mathcal{S}$ be a class of separable, unital $\mathrm{C}^{*}$-algebras which can be finitely presented with weakly stable relations, which is closed under taking direct sums and which contains all finite dimensional $\mathrm{C}^{*}$-algebras. Let $A$ be a separable, simple, unital, stably finite, exact $\mathrm{C}^{*}$-algebra.

Then, $A \otimes \mathcal{Q}$ is TA$\mathcal{S}$ if and only if the following holds: There is $\eta > 0$ such that, for any $\epsilon > 0$ and any finite subset $\mathcal{F} \subset A \otimes \mathcal{Q}$, there are a projection $p \in A \otimes \mathcal{Q}$ and a unital $\mathrm{C}^{*}$-subalgebra $B \subset p(A \otimes \mathcal{Q})p$, $B \in \mathcal{S}$, such that
\begin{enumerate}
\item $\|pb - bp\| < \epsilon$ for all $b \in \mathcal{F}$,
\item $\dist(pbp,B) < \epsilon$ for all $b \in \mathcal{F}$,
\item $\tau({p}) > \eta$ for all $\tau \in T(A \otimes \mathcal{Q})$. 
\end{enumerate} 
\end{nprop}
\en

\bn
We will mainly be interested in the cases where $\mathcal{S}$ is the class of finite dimensional $\mathrm{C}^{*}$-algebras or the class of interval algebras, i.e., in TAF and TAI algebras, respectively. For such algebras, if they are in addition nuclear and satisfy the Universal Coefficient Theorem (UCT), the Elliott invariant is complete. Moreover, we know the range of the invariant in these situations. We summarize the relevant results from \cite{Lin:TAFduke, Lin:asu-class} for the reader's convenience.

\begin{ntheorem}
Let $A_{i}$, $i=0,1$, be separable, simple, unital, nuclear $\mathrm{C}^{*}$-algebras which satisfy the UCT. Suppose $A_{1}$ and $A_{2}$ are TAI. 

Then, $A_{1} \cong A_{2}$ if and only if their Elliott invariants
\[
(\mathrm{K}_{0}(A_{i}), \mathrm{K}_{0}(A_{i})_{+}, [1_{A_{i}} ], \mathrm{K}_{1}(A_{i}), T(A_{i}), r_{A_{i}}: T(A_{i}) \to S(\mathrm{K}_{0}(A_{i})))
\]
are isomorphic, and every isomorphism of invariants lifts to a $^{*}$-isomorphism of algebras. 

Moreover, the $A_{i}$ are approximately subhomogeneous (ASH) algebras of topological dimension (and hence decomposition rank) at most 2, and are approximately homogeneous (AH)  of topological dimension at most 3. 

Finally, they are in fact TAF if and only if they have real rank zero, and in this case the classifying invariant degenerates to ordered $\mathrm{K}$-theory,
\[
(\mathrm{K}_{0}(A_{i}), \mathrm{K}_{0}(A_{i})_{+}, [1_{A_{i}} ], \mathrm{K}_{1}(A_{i})).
\]
\end{ntheorem}
\en

\bn
\label{H}
The previous theorem illustrates that classification in term of $\mathrm{K}$-theory works best in the case of real rank zero, i.e., with an abundance of projections around. In \cite{Win:localizingEC} it was shown that classification up to UHF-stability will yield classification up to $\mathcal{Z}$-stability, where $\mathcal{Z}$ denotes the Jiang--Su algebra, cf.\ \cite{JiaSu:Z, RorWin:Z-revisited}.  This is useful when there are at least enough projections to distinguish traces, since then tensoring with UHF algebras will enforce real rank zero. In this case it only remains to confirm $\mathcal{Z}$-stability, for which we have all sorts of highly useful criteria (for example finite decomposition rank). We summarize the situation as follows.   

\begin{ntheorem}
Let $A$ be a separable, simple, unital $\mathrm{C}^{*}$-algebra with finite decomposition rank. 

Then, conditions {\rm (i)}--{\rm (iv)} below are equivalent:
\begin{enumerate}
\item The canonical map $T(A) \to S(\mathrm{K}_{0}(A))$ is a homeomorphism.
\item $\mathrm{K}_{0}(A)$ separates the tracial states of $A$.
\item $A$ is rationally TAF, i.e., $A \otimes \mathcal{Q}$ is TAF.
\item $A$ is rationally of real rank zero, i.e., $A \otimes \mathcal{Q}$ has real rank zero. 
\end{enumerate}

If, moreover, $A$ satisfies the UCT, under any of these conditions $A$ is ASH of topological dimension at most 2, and
such $\mathrm{C}^{*}$-algebras are classified by their ordered $\mathrm{K}$-theory.
\end{ntheorem}
\en

\section{From approximate tracial embeddings to tracial approximations}
\label{section2}

\noindent
In this section we introduce and illustrate the `classification by embedding' method. We start with a technical result that allows us to compare order zero maps in terms of traces. Recall that a completely positive contractive (c.p.c.) map is order zero if it preserves orthogonality. There is a structure theorem for such maps which in particular yields a notion of functional calculus, see \cite{WinZac:order-zero} for a detailed exposition. In the sequel, we will encounter matrix algebras of different sizes. We will usually write matrix units for these in the form $e_{mn}$, without distinguishing between the matrix sizes; this should cause no confusion.

\bn
\label{E}
\begin{nprop}
Let $A$ be a separable, simple, unital $\mathrm{C}^{*}$-algebra with strict comparison. Let $F$ be a finite dimensional $\mathrm{C}^{*}$-algebra and let 
\[
\varphi:F \to A,
\]
\[
\varphi_{i}:F \to A, \; i \in \mathbb{N},
\]
be c.p.c.\ order zero maps such that, for each $x \in F_{+}$ and $f \in \mathcal{C}_{0}((0,1])_{+}$,
\[
\sup \{|\tau(f(\varphi)(x) - f(\varphi_{i})(x))| \mid \tau \in T(A) \} \stackrel{i \to \infty}{\longrightarrow} 0
\]
and 
\[
\limsup_{i} \|f(\varphi_{i})(x)\| \le \|f(\varphi)(x)\|.
\]

Then, there are
\[
s_{i} \in (M_{4} \otimes A)^{1}, i \in \mathbb{N},
\]
such that for each $y \in F_{+}$
\begin{equation}
\label{E3}
\|s_{i}(1_{4}\otimes \varphi(y)) - (e_{11} \otimes \varphi_{i}(y)) s_{i}\| \stackrel{i \to \infty}{\longrightarrow} 0
\end{equation}
and
\begin{equation}
\label{E4}
\|(e_{11}\otimes \varphi_{i}(y)) s_{i}s_{i}^{*} - (e_{11} \otimes \varphi_{i}(y)) \| \stackrel{i \to \infty}{\longrightarrow} 0,
\end{equation}
and that
\begin{equation}
\label{E1}
s_{i}s_{i}^{*} \in \overline{ (e_{11} \otimes \varphi_{i}(1_{F})) M_{4} \otimes A  (e_{11} \otimes \varphi_{i}(1_{F}))}
\end{equation}
and
\begin{equation}
\label{E2}
s_{i}^{*}s_{i} \in \overline{ (1_{4} \otimes \varphi(1_{F})) M_{4} \otimes A  (1_{4} \otimes \varphi(1_{F}))}.
\end{equation}
\end{nprop}

\begin{nproof}
We first show the statement in the case $F = \mathbb{C}$. When checking \eqref{E3} and \eqref{E4} it will clearly suffice to consider $y = 1_{\mathbb{C}}$; let us write $h$ and $h_{i}$ for $\varphi(1_{\mathbb{C}})$ and $\varphi_{i}(1_{\mathbb{C}})$, respectively.

For the moment let us fix $L \in \mathbb{N}$ and $\epsilon >0$. Define functions
\[
d^{({c})}_{l}, \tilde{d}^{({c})}_{l}, \hat{d}^{({c})}_{l} \in \mathcal{C}_{0}((0,1])^{1}_{+} \mbox{ for } c \in \{0,1\}, l \in \{1,\ldots,L\},
\] 
with the following properties:
\begin{itemize}
\item[(a)] $\|\tilde{d}^{({c})}_{l} \mathrm{id}_{(0,1]} - \frac{l}{L} \tilde{d}^{({c})}_{l}\| < \frac{1}{L}$, $c \in \{0,1\}$, $l \in \{1,\ldots, L\}$
\item[(b)] $d^{({c})}_{l} \tilde{d}^{({c})}_{l} = d^{({c})}_{l}$,  $c \in \{0,1\}$, $l \in \{1,\ldots, L\}$
\item[({c})] $\tilde{d}^{({c})}_{l} \tilde{d}^{({c})}_{l'} = 0$,  $c \in \{0,1\}$, $l \neq l' \in \{1,\ldots, L\}$
\item[(d)] $\hat{d}^{({c})}_{l} = (d^{({c})}_{l} - \epsilon)_{+}$, $c \in \{0,1\}$, $l \in \{1,\ldots, L\}$
\item[(e)] $\| \mathrm{id}_{(0,1]} - \sum_{c=0}^{1} \sum_{l=1}^{L} \frac{l}{L} \cdot d^{({c})}_{l}\| \le \frac{1}{L}$.
\end{itemize}
For $c \in \{0,1\}$, set
\[
N^{({c})}:= \{l \in \{1,\ldots,L\} \mid d_{l}^{({c})}(h) \neq 0\}
\]
and
\[
N_{\times}^{({c})} := \{1,\ldots,L\} \setminus N^{({c})}.
\]
One checks that there is $\bar{\imath} \in \mathbb{N}$ such that for each $i \ge \bar{\imath}$
\[
\tau(\tilde{d}^{({c})}_{l}(h_{i})) \le \frac{3}{2} \tau(\tilde{d}^{({c})}_{l}(h))
\]
for all $ \tau \in T(A)$, $c \in \{0,1\}$, $l \in N^{({c})}$ and 
\[
\|\tilde{d}^{({c})}_{l}(h_{i})\| < \epsilon
\]
for $c \in \{0,1\}$, $l \in N_{\times}^{({c})}$. But then
\[
d_{\tau}(d^{({c})}_{l}(h_{i})) \le \tau(\tilde{d}^{({c})}_{l}(h_{i})) \le \frac{3}{2} \tau(\tilde{d}^{({c})}_{l}(h)) \le \frac{3}{2} d_{\tau}(\tilde{d}^{({c})}_{l}(h))
\]
for $i \ge \bar{\imath}$, $\tau \in T(A)$, $c \in \{0,1\}$, $l \in N^{(c)}$.

Now by comparison,
\[
\langle d^{({c})}_{l}(h_{i})\rangle \le 2 \cdot \langle\tilde{d}^{({c})}_{l}(h) \rangle
\]
in $\mathrm{Cu}(A)$ for $i \ge \bar{\imath}$, $c \in \{0,1\}$, $l \in N^{({c})}$. 

By the Kirchberg--R{\o}rdam Lemma \cite[Lemma~2.2]{KirRor:pi2} there are
\[
s^{({c})}_{i,l} \in (M_{2} \otimes A)^{1}
\]
such that
\[
s^{({c})}_{i,l} s^{({c})*}_{i,l} = e_{11} \otimes \hat{d}^{({c})}_{l}(h_{i})
\]
and
\[
s^{({c})*}_{i,l} s^{({c})}_{i,l} \in \mathrm{her}(1_{2} \otimes \tilde{d}^{({c})}_{l}(h)) \subset M_{2} \otimes A
\]
for $ i \ge \bar{\imath}$, $c \in \{0,1\}$, $l \in N^{({c})}$. Set
\[
s^{(l)}_{i,l}:= 0
\]
for $i \ge \bar{\imath}$, $c \in \{0,1\}$, $l \in N_{\times}^{({c})}$.

For $ i \ge \bar{\imath}$, $c \in \{0,1\}$ define 
\[
\textstyle
s^{({c})}_{i}:= \sum_{l=1}^{L} s^{({c})}_{i,l} \in M_{2} \otimes A,
\]
then
\[
\textstyle
\| s^{({c})}_{i}s^{({c})*}_{i} - \sum_{l=1}^{L} e_{11} \otimes \hat{d}^{({c})}_{l}(h_{i}) \| \le \epsilon
\]
and
\[
\textstyle
\|(e_{11} \otimes h_{i}) s^{({c})}_{i} s^{({c})*}_{i} - \sum_{l=1}^{L} \frac{l}{L}\cdot e_{11} \otimes d^{({c})}_{l}(h_{i}) \| \le 2 \epsilon + 1/L,
\]
whence
\[
\textstyle
\|(e_{11} \otimes h_{i})(\sum_{c=0}^{1} s^{({c})}_{i}s^{({c})*}_{i}) - (e_{11} \otimes h_{i})\| \le 3 \epsilon + 3/L
\]
for $i \ge \bar{\imath}$.

Moreover, 
\begin{eqnarray*}
\lefteqn{ \| s^{({c})}_{i} (1_{2} \otimes h) - (e_{11} \otimes h_{i}) s^{({c})}_{i} \| }\\
& = & \| \textstyle \sum_{l=1}^{L} s^{({c})}_{i,l} (1_{2} \otimes h) - \sum_{l=1}^{L}  (e_{11} \otimes h_{i}) s^{({c})}_{i,l} \|\\
& \le & \| \textstyle \sum_{l=1}^{L} \frac{l}{L} \cdot s^{({c})}_{i,l}  - \sum_{l=1}^{L}  (e_{11} \otimes h_{i}) s^{({c})}_{i,l} \| + 1/L\\
& \le & 1/L + 2 \epsilon
\end{eqnarray*}
for $c \in \{0,1\}$, $i \ge \bar{\imath}$.

Now define
\[
\textstyle
\tilde{s}_{i}:= \sum_{c=0}^{1} e_{1,c+1} \otimes s^{({c})}_{i} \in M_{2} \otimes M_{2} \otimes A,
\]
then
\[
\textstyle
\tilde{s}_{i} \tilde{s}^{*}_{i} = e_{11} \otimes \sum_{c=0}^{1} s^{({c})}_{i} s^{({c})*}_{i} \le 1_{M_{2} \otimes M_{2} \otimes A}
\]
and
\[
\|(e_{11} \otimes e_{11} \otimes h_{i}) \tilde{s}_{i} \tilde{s}^{*}_{i} - (e_{11} \otimes e_{11} \otimes h_{i}) \| \le 3 \epsilon + 3/L
\]
for $i \ge \bar{\imath}$.

Next, we check
\begin{eqnarray*}
\lefteqn{\|\tilde{s}_{i} (1_{2} \otimes 1_{2} \otimes h) - (e_{11} \otimes e_{11} \otimes h_{i}) \tilde{s}_{i} \|} \\
& = & \textstyle \| \sum_{c=0}^{1} e_{1,c+1} \otimes (s^{({c})}_{i}(1_{2} \otimes h)) \\
&&\textstyle - \sum_{c=0}^{1} e_{1,c+1} \otimes ((e_{11} \otimes h_{i}) s^{({c})}_{i}) \| \\
& \le & 2/L +4 \epsilon
\end{eqnarray*}
for $i \ge \bar{\imath}$.

Now if we let $L$ go to infinity and $\epsilon$ to zero, a diagonal sequence argument yields
\[
s_{i} \in (M_{2} \otimes M_{2} \otimes A)^{1}, \, i \in \mathbb{N},
\]
satisfying
\[
\| s_{i} (1_{2} \otimes 1_{2} \otimes h) - (e_{11} \otimes e_{11} \otimes h_{i}) s_{i} \| \stackrel{i \to \infty}{\longrightarrow} 0
\]
and
\[
\| (e_{11} \otimes e_{11} \otimes h_{i}) s_{i}s_{i}^{*} - e_{11} \otimes e_{11} \otimes h_{i} \| \stackrel{i \to \infty}{\longrightarrow} 0.
\]
Upon replacing $s_{i}$ with $(e_{11} \otimes h_{i}^{\frac{1}{i+1}})s_{i} (1_{4} \otimes h^{\frac{1}{i+1}})$ if necessary, we get \eqref{E1} and \eqref{E2}. We have thus verified the proposition when $F= \mathbb{C}$.

\bigskip

Next suppose $F = M_{k}$ for some $k \in \mathbb{N}$. Let 
\[
\pi: M_{k} \longrightarrow A^{**}, \; \pi_{i}: M_{k} \longrightarrow A^{**}
\]
be supporting $^{*}$-homomorphisms for $\varphi$ and for the $\varphi_{i}$, respectively; cf.\ \cite{WinZac:order-zero}. Run the proposition for $\mathbb{C} \cong e_{11} M_{k} e_{11}$ and for $\varphi|_{e_{11} M_{k} e_{11}}$ and  $\varphi_{i}|_{e_{11} M_{k} e_{11}}$; denote the resulting elements of $(M_{4} \otimes A)^{1}$ by $\tilde{s}_{i}$. The `amplified' elements 
\[
s_{i}:= \textstyle \sum_{m=1}^{k} (e_{11} \otimes \pi_{i}(e_{m1})) \tilde{s}_{i} (1_{4} \otimes \pi(e_{1m})) \in (M_{4} \otimes A)^{1}
\]
will then satisfy \eqref{E3} through \eqref{E2} above. This verifies the proposition in the case $F = M_{k}$. 

\bigskip

When $F$ is a sum of $N$ matrix algebras, run the proposition for each matrix block separately to obtain elements $s_{i}^{(1)}, \ldots, s_{i}^{(N)} \in (M_{4} \otimes A)^{1}$. By \eqref{E1} and \eqref{E2}, the elements $s_{i}^{(1)}(s_{i}^{(1)})^{*}, \ldots , s_{i}^{(N)}(s_{i}^{(N)})^{*}$ are pairwise orthogonal for each fixed $i$, and the same goes for the $(s_{i}^{(1)})^{*}s_{i}^{(1)}, \ldots , (s_{i}^{(N)})^{*}s_{i}^{(N)}$. Therefore we can define elements
\[
s_{i}:= s_{i}^{(1)} + \ldots + s_{i}^{(N)} \in (M_{4} \otimes A)^{1}.
\]
These will again satisfy \eqref{E3} through \eqref{E2}, thus verifying the proposition for an arbitrary finite dimensional $\mathrm{C}^{*}$-algebra $F$. 
\end{nproof}
\en

\bn
\label{A}
\begin{ntheorem}
Let $\mathcal{S}$ be a class of separable, unital $\mathrm{C}^{*}$-algebras which can be finitely presented with weakly stable relations. Suppose further that $\mathcal{S}$ is closed under taking direct sums and under taking tensor products with finite dimensional $\mathrm{C}^{*}$-algebras, and that $\mathcal{S}$ contains all finite dimensional $\mathrm{C}^{*}$-algebras.

Let $A$ be a separable, simple, unital $\mathrm{C}^{*}$-algebra with $\dim_{\mathrm{nuc}} A < \infty$ and $T(A) \neq \emptyset$, and let
\[
\big(A \stackrel{\sigma_{i}}{\longrightarrow} B_{i}  \stackrel{\varrho_{i}}{\longrightarrow} A \big)_{i \in \mathbb{N}}
\]
be a system of maps with the following properties:

\begin{enumerate}
\item $B_{i} \in \mathcal{S}$, $i \in \mathbb{N}$
\item $\varrho_{i}$ is an embedding for each $i \in \mathbb{N}$
\item $\sigma_{i}$ is c.p.c.\ for each $i \in \mathbb{N}$
\item $\bar{\sigma}:A \to \prod_{\mathbb{N}} B_{i} / \bigoplus_{\mathbb{N}} B_{i}$ induced by the $\sigma_{i}$ is a unital $^{*}$-homomorphism
\item $\sup\{ |\tau(\varrho_{i}\sigma_{i}(a) - a)| \mid \tau \in T(A)\} \stackrel{i \to \infty}{\longrightarrow} 0$ for each $a \in A$.
\end{enumerate}

Then, $A \otimes \mathcal{Q}$ is TA$\mathcal{S}$.
\end{ntheorem}

\begin{nproof}
Let 
\[
(F_{j}= F_{j}^{(0)} \oplus \ldots \oplus F_{j}^{(m)}, \psi_{j},\varphi_{j})_{j \in \mathbb{N}}
\]
be a system of $m$-decomposable c.p.\ approximations for $A$ with c.p.c.\ approximately order zero maps $\psi_{j}$ and c.p.c.\ order zero maps $\varphi^{(l)}_{j}= \varphi_{j}|_{F^{(l)}_{j}}$ as in \cite[Proposition~4.2]{Winter:dimnuc-Z-stable}; see also \cite[Proposition~4.3]{WinterZac:dimnuc}. (Here, `approximately order zero' means that for $a,b \in A$ with $ab=0$ we have $\psi_{j}(a) \psi_{j}(b) \longrightarrow 0$ as $j$ goes to infinity.)

By weak stability of order zero maps, for each $j \in \mathbb{N}$ there are c.p.c.\ order zero maps 
\[
\tilde{\varphi}_{j,i}^{(l)}:F^{(l)}_{j} \to B_{i}, \, i \in \mathbb{N}, l \in \{0,\ldots,m\},
\]
such that
\[
\| \tilde{\varphi}^{(l)}_{j,i}(x) - \sigma_{i}\varphi^{(l)}_{j}(x) \| \stackrel{i\to \infty}{\longrightarrow} 0,\, x \in F^{(l)}_{j}.
\]
Now for each $x \in (F^{(l)}_{j})_{+}$ and $f \in \mathcal{C}_{0}((0,1])_{+}$, we have
\[
\| f(\tilde{\varphi}_{j,i}^{(l)})(x) - \sigma_{i}f(\varphi^{(l)}_{j})(x) \| \stackrel{ i \to \infty}{\longrightarrow} 0,
\]
hence
\[
\sup_{\tau \in T(A)} | \tau(f(\varrho_{i} \tilde{\varphi}_{j,i}^{(l)})(x) - f(\varphi_{j}^{(l)})(x)) | \stackrel{ i \to \infty}{\longrightarrow} 0.
\]
and
\[
\| \varrho_{i}f(\tilde{\varphi}_{j,i}^{(l)})(x) - \varrho_{i}\sigma_{i}f(\varphi^{(l)}_{j})(x) \| \stackrel{ i \to \infty}{\longrightarrow} 0.
\]
Since the $\varrho_{i}$ are embeddings, the $\sigma_{i}$ are (eventually) nonzero and approximately multiplicative, we see that 
\[
\limsup_{i} \|\varrho_{i}\sigma_{i} f(\varphi^{(l)}_{j})(x)\| = \|f(\varphi^{(l)}_{j})(x)\|,
\]
whence
\[
\limsup_{i} \|f(\varrho_{i} \tilde{\varphi}^{(l)}_{j,i})(x)\| \le \|f(\varphi^{(l)}_{j})(x)\|
\]
for $x \in (F^{(l)}_{j})_{+}$, $f \in \mathcal{C}_{0}((0,1])_{+}$.

By Proposition~\ref{E}, there are 
\[
s_{j,i}^{(l)} \in (M_{4} \otimes A)^{1}, \, i \in \mathbb{N},
\]
such that
\[
\| s^{(l)}_{j,i} (1_{4} \otimes \varphi^{(l)}_{j}(x)) - (e_{11} \otimes \varrho_{i} \tilde{\varphi}_{j,i}^{(l)}(x)) s_{j,i}^{(l)} \| \stackrel{ i \to \infty}{\longrightarrow} 0
\]
and
\[
\| (e_{11} \otimes \varrho_{i} \tilde{\varphi}_{j,i}^{(l)}(x)) s_{j,i}^{(l)}s_{j,i}^{(l)*} - e_{11} \otimes \varrho_{i} \tilde{\varphi}^{(l)}_{j,i}(x) \| \stackrel{ i \to \infty}{\longrightarrow} 0
\]
for each $x \in F^{(l)}_{j}$.

We obtain contractions
\[
s^{(l)}_{j} \in (M_{4} \otimes A)_{\infty} \cong M_{4} \otimes A_{\infty}
\]
with
\[
s^{(l)}_{j} (1_{4} \otimes \iota \varphi^{(l)}_{j}(x)) = (e_{11} \otimes \bar{\varrho}\bar{\sigma} \varphi^{(l)}_{j}(x)) s_{j}^{(l)}
\]
and
\[
(e_{11} \otimes \bar{\varrho} \bar{\sigma} \varphi^{(l)}_{j}(x)) s^{(l)}_{j}s^{(l)*}_{j} = e_{11} \otimes \bar{\varrho}\bar{\sigma} \varphi^{(l)}_{j}(x),
\]
where
\[
\textstyle
\bar{\varrho}: \prod B_{i}/\bigoplus B_{i} \to A_{\infty}
\]
is the $^{*}$-homomorphism induced by the $\varrho_{i}$ and
\[
\iota:A \to A_{\infty}
\]
is the canonical embedding. Let
\[
\bar{\iota}: A_{\infty} \to (A_{\infty})_{\infty}
\]
be induced by the canonical embedding
\[
\iota:A \to A_{\infty},
\]
i.e.,
\[
[(a_{j})_{j \in \mathbb{N}}] \stackrel{\bar{\iota}}{\mapsto} [(\iota(a_{j}))_{j \in \mathbb{N}}].
\]
Let 
\[
\bar{\gamma}: A_{\infty} \to (A_{\infty})_{\infty}
\]
be the $^{*}$-homomorphism induced by
\[
\bar{\varrho}\bar{\sigma}: A \to A_{\infty}.
\]
Let
\[
\textstyle
\bar{\varphi}^{(l)}: \prod_{j} F^{(l)}_{j}/\bigoplus_{j} F^{(l)}_{j} \to A_{\infty}
\]
and
\[
\textstyle
\bar{\psi}^{(l)}: A \to \prod_{j} F_{j}^{(l)}/\bigoplus_{j} F^{(l)}_{j}
\]
be the maps induced by the $\varphi_{j}^{(l)}$ and $\psi_{j}^{(l)}$, respectively; these will automatically be c.p.c.\ order zero.

Define
\[
\bar{s}^{(l)}:= [(s^{(l)}_{j})_{j \in \mathbb{N}} ] \in (M_{4} \otimes A_{\infty})_{\infty} \cong M_{4} \otimes (A_{\infty})_{\infty},
\]
then
\[
\bar{s}^{(l)} (1_{4} \otimes \bar{\iota} \bar{\varphi}^{(l)} \bar{\psi}^{(l)}(a)) = (e_{11} \otimes \bar{\gamma} \bar{\varphi}^{(l)} \bar{\psi}^{(l)}(a)) \bar{s}^{(l)}
\]
and
\[
(e_{11} \otimes \bar{\gamma} \bar{\varphi}^{(l)} \bar{\psi}^{(l)}(a)) \bar{s}^{(l)}\bar{s}^{(l)*} = e_{11} \otimes \bar{\gamma} \bar{\varphi}^{(l)} \bar{\psi}^{(l)}(a).
\]
Note that for each $a \in A$
\[
\bar{\varphi}^{(l)} \bar{\psi}^{(l)}(1_{A}) \iota(a) = \bar{\varphi}^{(l)} \bar{\psi}^{(l)}(a)
\]
by \cite[Proposition~4.2]{Winter:dimnuc-Z-stable}, and so in particular
\[
(\bar{\varphi}^{(l)} \bar{\psi}^{(l)}(1_{A}))^{\frac{1}{2}}\iota(a) \in \mathrm{C}^{*}(\bar{\varphi}^{(l)}\bar{\psi}^{(l)}(b) \mid b \in A),
\]
which in turn implies that
\begin{eqnarray*}
\lefteqn{ \bar{s}^{(l)} (1_{4} \otimes (\bar{\iota} \bar{\varphi}^{(l)} \bar{\psi}^{(l)} (1_{A}))^{\frac{1}{2}})  (1_{4} \otimes \bar{\iota} \iota (a))  }\\
& = & (e_{11} \otimes \bar{\gamma}(\iota(a)) (\bar{\gamma} \bar{\varphi}^{(l)}\bar{\psi}^{(l)}(1_{A}))^{\frac{1}{2}}) \bar{s}^{(l)} \\
& = & (e_{11} \otimes \bar{\gamma}(\iota(a))) (e_{11} \otimes (\bar{\gamma} \bar{\varphi}^{(l)}\bar{\psi}^{(l)}(1_{A}))^{\frac{1}{2}}) \bar{s}^{(l)}.
\end{eqnarray*}
Set 
\begin{eqnarray*}
\bar{v} & := & \textstyle \sum_{l=1}^{m+1} e_{1,l} \otimes ((e_{11} \otimes (\bar{\gamma} \bar{\varphi}^{(l)} \bar{\psi}^{(l)} (1_{A}))^{\frac{1}{2}})\bar{s}^{(l)})\\
& = &  \textstyle \sum_{l=1}^{m+1} e_{1,l} \otimes (\bar{s}^{(l)}    (1_{4} \otimes (\bar{\iota} \bar{\varphi}^{(l)} \bar{\psi}^{(l)} (1_{A}))^{\frac{1}{2}})) \in M_{m+1} \otimes M_{4} \otimes (A_{\infty})_{\infty},
\end{eqnarray*}
then
\[
\textstyle
\bar{v} \bar{v}^{*} = \sum_{l=1}^{m+1} e_{11} \otimes e_{11} \otimes \bar{\gamma}\bar{\varphi}^{(l)} \bar{\psi}^{(l)}(1_{A}) = e_{11} \otimes e_{11} \otimes \bar{\gamma}(1_{A_{\infty}}),
\]
so in particular $\bar{v}$ is a partial isometry.

Moreover, we check that for $a \in A$
\begin{eqnarray*}
\lefteqn{ \bar{v} (1_{m+1} \otimes 1_{4} \otimes \bar{\iota} \iota (a))} \\
 & = & \textstyle \sum_{l=1}^{m+1} e_{1,l} \otimes ((e_{11} \otimes \bar{\gamma}(\iota(a))) (e_{11} \otimes (\bar{\gamma} \bar{\varphi}^{(l)} \bar{\psi}^{(l)}(1_{A}))^{\frac{1}{2}})\bar{s}^{(l)}) \\
 & = & (e_{11} \otimes e_{11} \otimes \bar{\gamma}(\iota(a))) \bar{v},
\end{eqnarray*}
whence for $a \in A$
\begin{eqnarray*}
\bar{v}^{*} \bar{v} (1_{m+1} \otimes 1_{4} \otimes \bar{\iota} \iota(a)) & = & \bar{v}^{*} ( e_{11} \otimes e_{11} \otimes \bar{\gamma}(\iota(a))) \bar{v} \\
& = & (1_{m+1} \otimes 1_{4} \otimes \bar{\iota} \iota(a)) \bar{v}^{*} \bar{v}
\end{eqnarray*}
in $M_{m+1} \otimes M_{4} \otimes (A_{\infty})_{\infty}$.

Now for every finite subset $\mathcal{F} \subset A^{1}_{+}$ and $\epsilon>0$ there are $i \in \mathbb{N}$ and $v \in M_{m+1} \otimes M_{4} \otimes A$ such that
\begin{itemize}
\item[(a)] $vv^{*} = e_{11} \otimes e_{11} \otimes \varrho_{i}(1_{B_{i}})$
\item[(b)] $(\mathrm{tr}_{M_{m+1} \otimes M_{4}} \otimes \tau) (vv^{*}) \ge \frac{1}{2(m+1)4}$ for all $ \tau \in T(A)$
\item[({c})] $\| [v^{*}v, 1_{m+1} \otimes 1_{4} \otimes a] \| < \epsilon $ for all $ a \in \mathcal{F}$
\item[(d)] $ \| v^{*} v (1_{m+1} \otimes 1_{4} \otimes a) - v^{*} (e_{11} \otimes e_{11} \otimes \varrho_{i} \sigma_{i}(a))v\| < \epsilon$ for all $ a \in \mathcal{F}$.
\end{itemize}

Define
\[
\kappa: B_{i} \to M_{m+1} \otimes M_{4} \otimes A
\]
by
\[
\kappa(b) := v^{*}(e_{11} \otimes e_{11} \otimes \varrho_{i}(b)) v,
\]
then $\kappa$ is an embedding such that
\begin{itemize}
\item[(e)] $(\mathrm{tr}_{M_{m+1} \otimes M_{4}} \otimes \tau) (1_{\kappa(B_{i})}) \ge \frac{1}{2(m+1)4}$ for all $ \tau \in T(A)$
\item[(f)] $\| [1_{\kappa(B_{i})}, 1_{m+1} \otimes 1_{4} \otimes a] \| < \epsilon $ for all $ a \in \mathcal{F}$
\item[(g)] $ 1_{\kappa(B_{i})} (1_{m+1} \otimes 1_{4} \otimes a) 1_{\kappa(B_{i})} \in_{\epsilon} \kappa(B_{i})$ for all $ a \in \mathcal{F}$.
\end{itemize}
Now by Proposition~\ref{M}, $A \otimes \mathcal{Q}$ is TA$\mathcal{S}$.
\end{nproof}
\en

\bn
\label{N}
We note a Corollary which nicely illustrates the way in which we are going to use Theorem~\ref{A} towards classification results in the subsequent sections. Our method reduces the problem of showing that a $\mathrm{C}^{*}$-algebra $A$ is classifiable to the problem of embedding it into a classifiable (TAF, or TAI) $\mathrm{C}^{*}$-algebra $B$ such that the embedding induces an isomorphism at the level of invariants. Then, one needs to lift the inverse of this isomorphism to an embedding of the classifiable model into the original algebra. Theorem~\ref{A} now moves the image of the composition of these two embeddings into a position compatible with the TAF (or TAI) condition for $A$. 

There are quite a few tools available for finding such embeddings: In our applications, the map from $A$ to $B$ will usually come from properties related to quasidiagonality. The source of the map from $B$ to $A$ depends on the situation; for example, it might come from the existence theorem of \cite{Rob:NCCW}.   

Although in our applications we are not in the exact situation of the Corollary, it will serve as a blueprint for tracially approximate versions which essentially follow the same pattern, but which do not require keeping track of the entire $\mathrm{K}$-theoretic information of $A$ (see \ref{B}, \ref{C} and \ref{F} below).   

Recall from \cite[Section~3.1]{Rob:NCCW} that for a unital $\mathrm{C}^{*}$-algebra $D$ the ordered semigroup $\mathrm{Cu}^{\sim}(D)$ is given by formal differences $x - n \cdot [1_{D}]$, with $x \in \mathrm{Cu}(D)$ and $n \in \mathbb{N}$. More precisely, $\mathrm{Cu}^{\sim}(D)$ is the quotient of $\mathrm{Cu}(D) \times \mathbb{N}$ by the equivalence relation given by $(x,n) \sim (y,m) : \Leftrightarrow x + m \cdot [1_{D}] + k \cdot [1_{D}] = y + n \cdot [1_{D}] + k \cdot [1_{D}] $ in $\mathrm{Cu}(D)$ for some $k \in \mathbb{N}$.

\begin{ncor}
Let $A$ and $B$ be separable, simple, unital $\mathrm{C}^{*}$-algebras. Suppose that $\dim_{\mathrm{nuc}}A< \infty$ and that $B$ is an AI algebra. Suppose there is a unital embedding
\[
\sigma:A \to B
\]
inducing an isomorphism between the Elliott invariants
\[
(\mathrm{K}_{0}(A), \mathrm{K}_{0}(A)_{+}, [1_{A} ], \mathrm{K}_{1}(A), T(A), r_{A}: T(A) \to S(\mathrm{K}_{0}(A)))
\]
and
\[
(\mathrm{K}_{0}(B), \mathrm{K}_{0}(B)_{+}, [1_{B} ], \mathrm{K}_{1}(B), T(B), r_{B}: T(B) \to S(\mathrm{K}_{0}(B))).
\]

Then, $A \otimes \mathcal{Q}$ is TAI.
\end{ncor}

\begin{nproof}
From the K\"unneth Theorem it is clear that $\sigma \otimes \mathrm{id}_{\mathcal{Q}}$ also induces an isomorphism between the Elliott invariants of $A \otimes \mathcal{Q}$ and of $B \otimes \mathcal{Q}$. Under the conditions on $A \otimes \mathcal{Q}$ and $B \otimes \mathcal{Q}$ (they are both simple, unital, nuclear, $\mathcal{Z}$-stable, and have nonempty tracial state spaces, hence stable rank one), it follwos from \cite{AraPerToms:survey} that the Cuntz semigroups are determined in a natural way by the Elliott invariants, so that the isomorphim between the latter induces one between the former. But since also the classes of the units are preserved by this isomorphism, it then also induces an isomorphism $\mathrm{Cu}^{\sim}(\sigma \otimes \mathrm{id}_{\mathcal{Q}})$ between $\mathrm{Cu}^{\sim}(A \otimes \mathcal{Q})$ and $\mathrm{Cu}^{\sim}(B \otimes \mathcal{Q})$.  Since $B \otimes \mathcal{Q}$ is AI and $A \otimes \mathcal{Q}$ has stable rank one by \cite{Ror:Z-absorbing}, it follows from  \cite[Theorem~1]{Rob:NCCW}  that the inverse of this isomorphism lifts to a unital embedding 
\[
\varrho: B \otimes \mathcal{Q} \to A \otimes \mathcal{Q}.
\]
Let
\[
\big( B \otimes \mathcal{Q} \stackrel{\psi_{i}}{\longrightarrow} B_{i} \stackrel{\varphi_{i}}{\longrightarrow} B \otimes \mathcal{Q} \big)_{i \in \mathbb{N}}
\]
be a system of maps such that the $B_{i}$ are interval algebras, the $\psi_{i}$ are approximately multiplicative, the $\varphi_{i}$ are embeddings and $\varphi_{i} \psi_{i} \to \mathrm{id}_{B \otimes \mathcal{Q}}$ in point norm topology.  

For $i \in \mathbb{N}$ define 
\[
\sigma_{i}: A \otimes \mathcal{Q} \to B_{i} 
\]
and
\[
\varrho_{i}: B_{i} \to A \otimes \mathcal{Q}
\]
by 
\[
\sigma_{i} := \psi_{i} \circ (\sigma \otimes \mathrm{id}_{\mathcal{Q}}) 
\]
and
\[
\varrho_{i}:= \varrho \circ \varphi_{i}.
\]
It is now straightforward to check 
that the system 
\[
\big( A \otimes \mathcal{Q} \stackrel{\sigma_{i}}{\longrightarrow} B_{i} \stackrel{\varrho_{i}}{\longrightarrow} A \otimes \mathcal{Q} \big)_{i \in \mathbb{N}}
\]
satisfies the hypotheses of Theorem~\ref{A}, whence $A \otimes \mathcal{Q} \otimes \mathcal{Q} \cong A \otimes \mathcal{Q}$ is TAI.
\end{nproof}
\en

\bn
\label{B}
As a first incidence of a tracially approximate version of \ref{N}, in the monotracial situation we rediscover a special case of a striking recent result of Matui and Sato, \cite{MatSat:dr-UHF}. Our version is less general since we need to assume finite nuclear dimension; on the other hand, the proof is substantially simpler since it avoids the von Neumann algebra techniques pivotal for \cite{MatSat:dr-UHF} (see also \cite{MatuiSato:Comp}). We will see in the subsequent sections that our approach has the additional advantage that it applies in the situation of more general trace spaces.
 
\begin{ncor}
Let $A$ be a separable, simple, unital, monotracial $\mathrm{C}^{*}$-algebra with $\dim_{\mathrm{nuc}}A < \infty$. Suppose that $A$ is quasidiagonal.

Then, $A \otimes \mathcal{Q}$ is TAF.
\end{ncor}

\begin{nproof}
By \cite{BlaKir:limits}, $A \otimes \mathcal{Q}$ is NF, so there is a unital embedding
\[
\textstyle
\bar{\sigma}: A \otimes \mathcal{Q} \to \prod_{i} M_{n_{i}} / \bigoplus_{i} M_{n_{i}}
\]
for a suitable sequence $(n_{i})_{i \in \mathbb{N}} \subset \mathbb{N}$. Since $A \otimes \mathcal{Q}$ is nuclear, there is a sequence of c.p.c.\ maps
\[
\sigma_{i}: A \otimes \mathcal{Q} \to M_{n_{i}}, \; i \in \mathbb{N},
\] 
such that
\[
\tilde{\sigma}: A \otimes \mathcal{Q} \to \prod_{i} M_{n_{i}},
\]
given by
\[
\tilde{\sigma}(a) = (\sigma_{i}(a))_{i \in \mathbb{N}},
\]
lifts $\bar{\sigma}$.

Let
\[
\rho_{i}: M_{n_{i}} \to 1_{A} \otimes \mathcal{Q} \subset A \otimes \mathcal{Q}
\]
be a unital embedding, then clearly
\[
\tau(\rho_{i} \sigma_{i}(a) - a) \stackrel{i \to \infty}{\longrightarrow} 0
\]
for each $a \in A \otimes \mathcal{Q}$, where $\tau$ denotes the unique trace on $A \otimes \mathcal{Q}$.

Now Theorem~\ref{A} yields that $A \otimes \mathcal{Q} \otimes \mathcal{Q} \cong A \otimes \mathcal{Q}$ is TAF.
\end{nproof}
\en

\section{Free, minimal $\mathbb{Z}^{d}$-actions}
\label{section3}

\noindent
We now combine the method of the previous section with results of Lin and of Szab\'o to obtain our classification result for crossed products by $\mathbb{Z}^{d}$-actions.

\bn
\label{C}
\begin{ntheorem}
Let $A$ be a separable, simple, unital $\mathrm{C}^{*}$-algebra with $\dim_{\mathrm{nuc}}A< \infty$ and such that $A \otimes \mathcal{Q}$ has real rank zero; suppose the extreme boundary of the tracial state space of $A$, $\partial_{e}T(A)$, is nonempty and compact. 

Suppose further that for each $\tau \in \partial_{e}T(A)$ there are a simple, unital, monotracial AF algebra $D$ with trace $\delta$ and a unital embedding
\[
\alpha: A \to D
\]
with
\[
\delta \circ \alpha = \tau.
\]

Then, $A \otimes \mathcal{Q}$ is TAF.
\end{ntheorem}

\begin{nproof}
We may clearly replace $A$ by $A \otimes \mathcal{Q}$. 

For the moment fix $\mathcal{F} \subset (A \otimes \mathcal{Q})_{+}^{1}$ finite and $\epsilon >0$. 

Choose a finite partition of unity 
\[
(h_{\lambda})_{\Lambda}
\] 
for $\partial_{e}T(A)$ such that for each $\lambda \in \Lambda$ there is $\tau_{\lambda} \in \supp (h_{\lambda})$ such that 
\[
| \tau \otimes \tau_{\mathcal{Q}}(a) - \tau_{\lambda} \otimes \tau_{\mathcal{Q}}(a)| \le \epsilon
\]
for $\tau \in \supp (h_{\lambda})$, $a \in \mathcal{F}$.

Choose
\[
0 < \eta < \frac{\epsilon}{|\Lambda|}.
\]

Since $A$ has real rank zero, the image of the set of projections of $A$ in $\mathcal{C}(\partial_{e}T(A))_{+}^{1}$ is dense; using comparison one easily checks that there are pairwise orthogonal projections
\[
p_{\lambda} \in A, \, \lambda \in \Lambda,
\]
such that 
\[
|\tau(p_{\lambda}) - h_{\lambda}(\tau)| < \eta
\]
for $\lambda \in \Lambda$, $\tau \in \partial_{e}T(A)$. 

For each $\lambda \in \Lambda$, find $D_{\lambda}$, $\delta_{\lambda}$, $\alpha_{\lambda}$ as in the hypotheses, i.e., each $D_{\lambda}$ is simple, unital, AF with unique trace $\delta_{\lambda}$, and $\alpha_{\lambda}: A \otimes \mathcal{Q} \to D_{\lambda}$ is a unital embedding with $\delta_{\lambda} \circ \alpha_{\lambda} = \tau_{\lambda} \otimes \tau_{\mathcal{Q}}$, $\lambda \in \Lambda$. 

Find matrix algebras $M_{r_{\lambda,i}}, i \in \mathbb{N}$, and approximately multiplicative u.c.p.\ maps 
\[
\beta_{\lambda,i}: D_{\lambda} \to M_{r_{\lambda,i}};
\]
choose unital embeddings 
\[
\gamma_{\lambda,i}:M_{r_{\lambda,i}} \to \mathcal{Q}.
\]
Note that for each $a \in A \otimes \mathcal{Q}$
\[
\tau_{\mathcal{Q}} \circ \gamma_{\lambda,i} \circ \beta_{\lambda,i} \circ \alpha_{\lambda} (a) \stackrel{i \to \infty}{\longrightarrow} \tau_{\lambda} \otimes \tau_{\mathcal{Q}}(a).
\]

Next define
\[
\textstyle
\bar{B}_{i}:= \bigoplus_{\lambda \in \Lambda} M_{r_{\lambda,i}},
\]
\[
\bar{\sigma}_{i}: A \otimes \mathcal{Q} \to \bar{B}_{i},\, \bar{\sigma}_{i}:= \oplus_{\lambda} \beta_{\lambda,i} \circ \alpha_{\lambda},
\]
\[
\bar{\rho}_{i}: \bar{B}_{i} \to A \otimes \mathcal{Q}, \, \bar{\rho}_{i}:= \oplus_{\lambda} p_{\lambda} \otimes \gamma_{\lambda,i}.
\]

We check for $ a \in \mathcal{F}$ and $\tau \in \partial_{e}T(A)$
\begin{eqnarray*}
\lefteqn{ | \tau \otimes \tau_{\mathcal{Q}} (\bar{\rho}_{i} \bar{\sigma}_{i}(a) - a) |    } \\
& = &\textstyle |\sum_{\lambda} \tau(p_{\lambda}) \cdot \tau_{\mathcal{Q}} \circ \gamma_{\lambda,i} \circ \beta_{\lambda,i} \circ \alpha_{\lambda} (a) - h_{\lambda} (\tau) \cdot (\tau \otimes \tau_{\mathcal{Q}})(a) |\\
& \stackrel{i \to \infty}{\longrightarrow}& \textstyle |\sum_{\lambda} \tau(p_{\lambda}) \cdot (\tau_{\lambda} \otimes \tau_{\mathcal{Q}})  (a) - h_{\lambda} (\tau) \cdot (\tau \otimes \tau_{\mathcal{Q}})(a) |\\
& \le & \textstyle \sum_{\lambda} h_{\lambda}(\tau) \cdot  | (\tau_{\lambda} \otimes \tau_{\mathcal{Q}})  (a) - (\tau \otimes \tau_{\mathcal{Q}})(a) | + |\Lambda| \cdot \eta\\
& \le & \textstyle \sum_{\lambda} h_{\lambda}(\tau) \cdot \epsilon + |\Lambda| \cdot \eta \\
& \le & 2 \cdot \epsilon.
\end{eqnarray*}

Making $\mathcal{F}$ bigger and $\epsilon$ smaller will now produce 
\[
(A \otimes \mathcal{Q} \stackrel{\sigma_{i}}{\longrightarrow} B_{i} \stackrel{\rho_{i}}{\longrightarrow} A \otimes \mathcal{Q})_{i \in \mathbb{N}}
\]
as required for Theorem~\ref{A}.
\end{nproof}
\en

\bn
\label{D}
\begin{ncor}
Let $X$ be a compact metrizable space with finite covering dimension and $\beta: \mathbb{Z}^{d} \to \mathrm{Homeo}(X)$ a free, minimal action with compact space of ergodic measures.

Suppose that $\mathrm{K}_{0}(\mathcal{C}(X) \rtimes \mathbb{Z}^{d})$ separates the traces of $\mathcal{C}(X) \rtimes \mathbb{Z}^{d}$ (this is automatically satisfied if $X$ is a Cantor set).

Then, $(\mathcal{C}(X) \rtimes \mathbb{Z}^{d}) \otimes \mathcal{Q}$ is TAF. As a consequence, the crossed products themselves are classified by their ordered $\mathrm{K}$-theory.
\end{ncor}

\begin{nproof}
In \cite{Sza:dimnuc}, Szab\'o shows that $\mathrm{dim}_{\mathrm{nuc}}(\mathcal{C}(X) \rtimes \mathbb{Z}^{d}) < \infty$;  $\mathcal{C}(X) \rtimes \mathbb{Z}^{d}$ is simple since the action is minimal, hence $\mathcal{Z}$-stable by \cite[Corollary~6.3]{Winter:dimnuc-Z-stable}. Also, traces are separated by projections, and therefore $(\mathcal{C}(X) \rtimes \mathbb{Z}^{d}) \otimes \mathcal{Q}$ has real rank zero, see \cite{Ror:Z-absorbing}. 

By \cite[Theorem~9.3 and its proof]{Lin:crossed-product-AF-embedding}, for each $\tau \in T(\mathcal{C}(X) \rtimes \mathbb{Z}^{d})$ there is a unital embedding $\alpha: \mathcal{C}(X) \rtimes \mathbb{Z}^{d} \to D$ into a simple, unital, AF algebra with unique tracial state $\delta$ and such that $\tau = \delta \circ \alpha$. Now $(\mathcal{C}(X) \rtimes \mathbb{Z}^{d}) \otimes \mathcal{Q}$ is TAF by Theorem~\ref{C}. 

By \cite{LinNiu:KKlifting}, $\mathcal{C}(X) \rtimes \mathbb{Z}^{d}$ is classifiable ($\mathcal{C}(X) \rtimes \mathbb{Z}^{d}$ is well known to satisfy the UCT, see \cite{Bla:k-theory}); since it has real rank zero it is TAF.
\end{nproof}
\en

\section{Odd spheres}
\label{section4}

\noindent
Below we will see that $\mathrm{C}^{*}$-algebras associated to certain minimal homeomorphisms of spheres (of odd dimension at least 3) are classified by their spaces of invariant Borel probability measures; the crucial point here is that for the crossed products projections do not separate tracial states. The argument combines classification by embeddings with recent results of Strung and of Robert.

\bn
We start with a technical lemma, the crucial step of which relies on \cite[Lemma~4.4]{TikWin:Z-dr}. 
The result is implicitly contained in \cite[Sections~2.4 and 2.5]{Str:thesis}; a full proof is given below for the convenience of the reader. Recall that a \emph{purely} positive element of a $\mathrm{C}^{*}$-algebra is one that is not Cuntz equivalent to a projection.

\label{P}
\begin{nlemma}
Let $B= \mathcal{C}([0,1]) \otimes M_{r}$, $b \in B^{1}_{+}$, $\epsilon>0$ be given.

Then, there are $q, L \in \mathbb{N}$, purely positive elements $b_{1},\ldots,b_{L} \in (B \otimes M_{q})_{+}^{1}$ and numbers $\nu_{1},\ldots,\nu_{L} \in [0,1]$ such that the $b_{l}$ are pairwise orthogonal and such that 
\begin{equation}
\label{P1}
\textstyle
\big| \tau(b) - \sum_{l=1}^{L} \nu_{l} \cdot d_{\tau \otimes \mathrm{tr}_{q}}(b_{l})\big| < \epsilon
\end{equation}
for any $\tau \in T(B)$.
\end{nlemma}

\begin{nproof}
Once the $b_{l}$ are constructed, it will be enough to confirm \eqref{P1} for extremal traces of $B$, i.e.\ for traces of the form $\tau_{t} = \mathrm{ev}_{t} \otimes \mathrm{tr}_{r}$, $t \in [0,1]$. The numbers $d_{\tau_{t} \otimes \mathrm{tr}_{q}}(b_{l})$ are then just the ranks, divided by $r \cdot q$, of the matrices $b_{l}(t)$, $t \in [0,1]$. 

Moreover, $t \mapsto \tau_{t}(b)$ is just a positive continuous function of norm at most $1$ on $[0,1]$, and there are $2 \le L \in \mathbb{N}$, $0= t_{0}< t_{1}< \ldots< t_{L}=1$ and $\nu_{0},\ldots,\nu_{L} \in [0,1]$ such that for $t \in [t_{l-1},t_{l}]$, $l=1,\ldots,L$, we have  $|\tau_{t}(b) - \nu_{l}| < \epsilon/4$ and $|\nu_{l} - \nu_{l-1}|<\epsilon/4$. Set $t'_{l}:= (t_{l}-t_{l-1})/2$, $l=1,\ldots,L$.

Choose $q \in \mathbb{N}$ such that $1/q<\epsilon/4$. Let $D_{q} \subset M_{q}$ denote the subalgebra of diagonal matrices. For each $l=1,\ldots,L-1$, choose pairwise disjoint nondegenerate closed intervals $I_{l,1},\ldots, I_{l,q} \subset (t'_{l},t'_{l+1})$ and a function $a_{\frac{1}{2},l} \in \mathcal{C}([t'_{l},t'_{l+1}],D_{q})_{+}^{1}$ such $a_{\frac{1}{2},l}(t)$ has rank at most $1$ for each $t \in [t'_{l},t'_{l+1}]$ and such that for $t \in I_{l,s}$, the $s$-th diagonal entry of $a_{\frac{1}{2},l}(t)$ is $1$. 

Now by \cite[Lemma~4.4]{TikWin:Z-dr}, for each $l=1,\ldots,L-1$, there are 
\[
a_{0,l}, a_{1,l} \in \mathcal{C}([t'_{l},t'_{l+1}],D_{q})_{+}^{1}
\]
such that 
\[
a_{0,l} \perp a_{1,l}, \, a_{0,l} + a_{\frac{1}{2},l} + a_{1,l} = 1_{[t'_{l},t'_{l+1}]},  \mbox{ and } a_{0,l}(t'_{l}) = a_{1,l}(t'_{l+1}) = 1_{q}.
\]
Note that for each $t \in [t'_{l},t'_{l+1}]$, the ranks of $a_{0,l}(t)$ and $a_{1,l}(t)$ add up to at least $q-1$, and that $a_{0,l}(t'_{l+1}) = a_{1,l}(t'_{l}) = 0$.

We are now ready to define the $M_{r} \otimes M_{q}$-valued functions $b_{l}$, $l=1,\ldots,l$ as follows: For $l=2,\ldots,L-1$, set
\[
b_{l}(t):= 
\left\{
\begin{array}{ll}
1_{r} \otimes a_{1,l-1}(t), & t \in [t'_{l-1},t'_{l}] \\
1_{r} \otimes a_{0,l}(t), & t \in [t'_{l},t'_{l+1}] \\
0, & t \in [0,1] \setminus [t'_{l},t'_{l+1}].
\end{array}
\right.
\]
For $l=1,L$, set
\[
b_{1}(t):= 
\left\{
\begin{array}{ll}
1_{r} \otimes 1_{q}, & t \in [0,t'_{1}] \\
1_{r} \otimes a_{0,l}(t), & t \in [t'_{1},t'_{2}] \\
0, & t \in  [t'_{2},1]
\end{array}
\right.
\]
and
\[
b_{L}(t):= 
\left\{
\begin{array}{ll}
1_{r} \otimes 1_{q}, & t \in [t'_{L},1] \\
1_{r} \otimes a_{1,L-1}(t), & t \in [t'_{L-1},t'_{L}] \\
0, & t \in  [0,t'_{L-1}].
\end{array}
\right.
\]

It is clear from our construction that the $b_{l}$ are indeed continuous (thus well-defined) positive contractions; they each take the value $1_{r} \otimes 1_{q}$ for some $t$ and are $0$ for some other $t$, hence (the interval is connected) must be purely positive. Each $b_{l}$ is supported on $[t'_{l-1},t'_{l+1}]$. Moreover, for each $t \in [0,1]$, $b_{l}(t)$ can be nonzero for at most two values of $l$, which then must be consecutive. The $b_{l}$ are pairwise orthogonal since the $a_{0,l}$ and $a_{1,l}$ are orthogonal for each $l$.

Now let $t \in [t'_{\bar{l}},t'_{\bar{l}+1}]$ for some $\bar{l} \in \{1,\ldots,L-1\}$, so that $b_{l}(t) = 0$ for all $l \neq \bar{l},\bar{l}+1$. We estimate
\begin{eqnarray*}
\lefteqn{\textstyle | \tau_{t}(b) - \sum_{l=1}^{L} \nu_{l} \cdot d_{\tau_{t} \otimes \mathrm{tr}_{q}} (b_{l}) |} \\
& = & | \tau_{t}(b) -  \nu_{\bar{l}} \cdot d_{\tau_{t} \otimes \mathrm{tr}_{q}} (b_{\bar{l}})  -  \nu_{\bar{l}+1} \cdot d_{\tau_{t} \otimes \mathrm{tr}_{q}} (b_{\bar{l}+1}) | \\
& < & | \tau_{t}(b) -  \frac{\nu_{\bar{l}}}{r \cdot q} \cdot (\mathrm{rank} (b_{\bar{l}}(t))  +  \mathrm{rank} (b_{\bar{l}+1}(t))) | + \epsilon/4\\
& < & |\tau_{t}(b) - \nu_{\bar{l}} \cdot 1| + \epsilon/4 + 1/(r \cdot q) \\
& < & \epsilon.
\end{eqnarray*} 
For $t \in [0,t'_{1}]$ and $t \in [t'_{L},1]$ similar (in fact, easier) estimates hold, so we have confirmed \eqref{P1} for extremal traces of $B$.
\end{nproof}
\en

\bn
\label{F}
\begin{ntheorem}
Let $A$ and $B$ be separable, simple, unital $\mathrm{C}^{*}$-algebras. Suppose that $\dim_{\mathrm{nuc}}A< \infty$. Let $B$ be TAI and suppose there is a unital embedding
\[
\iota:A \to B
\]
such that
\[
T(\iota): T(B) \stackrel{\approx}{\longrightarrow} T(A)
\]
and such that
\[
\tau_{*} = \tau'_{*} \in S(\mathrm{K}_{0}(B)) \mbox{ for } \tau, \tau' \in T(B).
\]

Then, $A \otimes \mathcal{Q}$ is TAI.
\end{ntheorem}

\begin{nproof}
We may assume $A = A \otimes \mathcal{Q}$. Let $B_{i} \subset B$, $i \in \mathbb{N}$, be a sequence of interval algebras tracially approximating $B$ via u.c.p.\ maps 
\[
\psi_{i}:B \to B_{i},
\]
i.e., the $\psi_{i}$ are approximately multiplicative,
\[
\|\psi_{i}(b) - 1_{B_{i}}b1_{B_{i}}\| \stackrel{i \to \infty}{\longrightarrow} 0
\]
and
\[
\|[1_{B_{i}},b]\| \stackrel{i \to \infty}{\longrightarrow} 0,
\]
and
\[
\textstyle
\sup_{\tau \in T(B)} \{1 - \tau(1_{B_{i}}) \} < \epsilon_{i} 
\]
for some sequence $(\epsilon_{i})_{i \in \mathbb{N}} \subset (0,1)$, with $\epsilon_{i} \stackrel{i \to \infty}{\longrightarrow}0$.

Each $B_{i}$ is of the form
\[
\textstyle
B_{i} = \bigoplus_{j=1}^{M_{i}} B_{i,j}
\]
with each $B_{i,j}$ being nonzero and either $M_{r_{i,j}}$ or $\mathcal{C}([0,1]) \otimes M_{r_{i,j}}$ for some $r_{i,j}$. Set
\[
\lambda_{i,j}:= \tau(1_{B_{i,j}}) \in \mathbb{R}^{*}_{+}
\]
for some, hence all, $\tau \in T(B)$ (the $\lambda_{i,j}$ are nonzero by simplicity of $B$).

Choose $\mu_{i,j} \in \mathbb{Q}_{+}$ such that 
\[
\sup_{j \in \{1,\ldots,M_{i}\}} \{|1-{\textstyle \frac{\mu_{i,j}}{\lambda_{i,j}}}|\} < \frac{\epsilon_{i}}{3 \cdot M_{i} \cdot r_{i,j}}
\]
for $i \in \mathbb{N}$.

Let
\[
\gamma:\mathrm{Cu}(B \otimes \mathcal{Q}) \stackrel{\cong}{\longrightarrow} V(B\otimes \mathcal{Q}) \sqcup \mathrm{LAff}(T(B \otimes \mathcal{Q}))^{++}
\]
be the semigroup isomorphism (in $\mathrm{Mor}({\bf Cu})$) of \cite[Theorem~5.27]{AraPerToms:survey}. (The domain of the isomorphism in \cite{AraPerToms:survey} is $\mathrm{W}(B \otimes \mathcal{Q} \otimes \mathcal{K})$, which by \cite[Corollary~4.31 and Theorem~4.33]{AraPerToms:survey} can be identified with $\mathrm{Cu}(B \otimes \mathcal{Q})$.)

Let 
\[
\delta: \mathrm{LAff}(T(A \otimes \mathcal{Q}))^{++} \stackrel{\cong}{\longrightarrow} \mathrm{LAff}(T(B \otimes \mathcal{Q}))^{++}
\]
be induced by $\iota \otimes \mathrm{id}_{\mathcal{Q}}$, i.e.,
\[
\delta{f}(\tau) = f(\tau \circ (\iota \otimes \mathrm{id}_{\mathcal{Q}}))
\]
for $f \in \mathrm{LAff}(T(A \otimes \mathcal{Q}))^{++}$ and $\tau \in T(B \otimes \mathcal{Q})$. Note that 
\[
\delta^{-1}(g)(\tau \otimes \tau_{\mathcal{Q}}) = g(T(\iota)^{-1}(\tau) \otimes \tau_{\mathcal{Q}})
\]
for $g \in \mathrm{LAff}(T(B \otimes \mathcal{Q}))^{++}$ and $\tau \in T(A)$.

Let 
\[
\zeta: \mathbb{Q}_{+} \to V(\mathcal{Q}) \to V(A \otimes \mathcal{Q})
\] 
be the canonical map.

Define
\[
\kappa_{i}: \mathrm{Cu}(B_{i}) \to \mathrm{Cu}(A \otimes \mathcal{Q})
\]
by
\[
\textstyle
\kappa_{i}([b]) := \begin{cases}
\zeta(\frac{\mu_{i,j}}{\lambda_{i,j}} \cdot d_{\tau}([b])), &  \mbox{if } [b] \mbox{ is the class of a projection in } \mathcal{K} \otimes B_{i,j},\\
\delta^{-1}(\frac{\mu_{i,j}}{\lambda_{i,j}} \cdot \gamma([b \otimes 1_{\mathcal{Q}}])) ,& \mbox{if } b \mbox{ is purely positive in } \mathcal{K} \otimes B_{i,j};
\end{cases} 
\]
here, the dimension function $d_{\tau}$ comes from some $\tau \in T(B)$ (note that by our hypotheses on the pairing between $T(B)$ and $S(\mathrm{K}_{0}(B))$ $d_{\tau}([b])$ is independent of the particular choice of $\tau$ as long as $[b]$ is represented by a projection).

Note that $\frac{\mu_{i,j}}{\lambda_{i,j}} \cdot d_{\tau}([b]) \in \mathbb{Q}_{+}$ if $[b]$ is the class of a projection, hence $\kappa_{i}$ is well-defined. 

One checks that $\kappa_{i}$ is a semigroup homomorphism in $\mathrm{Mor}({\bf Cu})$. Let
\[
\kappa_{i}^{\sim}:\mathrm{Cu}^{\sim}(B_{i}) \to \mathrm{Cu}^{\sim}(A\otimes \mathcal{Q})
\]
be the induced map (cf.\ \cite[Section~3]{Rob:NCCW}), i.e.,
\[
\kappa_{i}^{\sim}([x] - n\cdot[1]) = \kappa_{i}([x]) - n\cdot [1];
\]
then $\kappa_{i}^{\sim} \in \mathrm{Mor}({\bf Cu})$. By \cite[Theorem~1]{Rob:NCCW}, $\kappa_{i}^{\sim}$ lifts to a unital $^{*}$-homomorphism
\[
\beta_{i}:B_{i} \to A \otimes \mathcal{Q}.
\]

If $B_{i,j}$ is a matrix algebra, then for a projection $p \in B_{i,j}$ and for traces $\tau \in T(A)$, $\tau' \in T(B)$,
\begin{eqnarray*}
\lefteqn{ | (\tau \otimes \tau_{\mathcal{Q}} )(\beta_{i}({p})) - (T(\iota)^{-1}(\tau)) ({p}) |  } \\
& = & | d_{\tau \otimes \tau_{\mathcal{Q}}}( \beta_{i}({p})) - \tau'({p}) | \\
& = & | d_{\tau \otimes \tau_{\mathcal{Q}}}( \kappa_{i}^{\sim}([p])) - d_{\tau'}([p]) | \\
& = & | d_{\tau \otimes \tau_{\mathcal{Q}}}( \kappa_{i}([p])) - d_{\tau'}([p]) | \\
& = & \textstyle | d_{\tau \otimes \tau_{\mathcal{Q}}}( \zeta(\frac{\mu_{i,j}}{\lambda_{i,j}}\cdot d_{\tau'}([p]))) - d_{\tau'}([p]) | \\
& = & \textstyle | \frac{\mu_{i,j}}{\lambda_{i,j}}\cdot d_{\tau'}([p]) - d_{\tau'}([p]) | \\
& \le &  \frac{\epsilon_{i}}{M_{i} \cdot r_{i,j}}.
\end{eqnarray*}
But then for each $b \in (B_{i,j})^{1}_{+}$ we have 
\[
| (\tau \otimes \tau_{\mathcal{Q}}) ( \beta_{i}(b)) - (T(\iota)^{-1}(\tau))(b) | < \frac{\epsilon_{i}}{M_{i}}
\]
for any $\tau \in T(A)$. 

Let us now consider the case where $B_{i,j}$ is of the form $\mathcal{C}([0,1]) \otimes M_{r_{i,j}}$.

We first check that $\beta|_{B_{i,j}}$ is injective, hence in particular maps purely positive elements to purely positive elements:

If $0 \neq b \in (B_{i,j})^{1}_{+}$ is purely positive, then for any $\tau \in T(B \otimes \mathcal{Q})$, 
\begin{eqnarray*}
\delta([\beta_{i}(b)])(\tau) & = & \delta(\kappa_{i}([b]))(\tau) \\
& = & \textstyle \delta \delta^{-1}(\frac{\mu_{i,j}}{\lambda_{i,j}} \cdot \gamma([b \otimes 1_{\mathcal{Q}}]))(\tau) \\
& = & \textstyle \frac{\mu_{i,j}}{\lambda_{i,j}} \cdot \gamma([b \otimes 1_{\mathcal{Q}}])(\tau) \\
& = & \textstyle \frac{\mu_{i,j}}{\lambda_{i,j}} \cdot d_{\tau}(b \otimes 1_{\mathcal{Q}}) \\
& \neq & 0,
\end{eqnarray*}
from which follows that $\beta|_{B_{i,j}}$ is injective.

Now let $b \in (B_{i,j})^{1}_{+}$ be arbitrary. With the aid of Lemma~\ref{P}, find $q, L \in \mathbb{N}$ and $b_{1},\ldots,b_{L} \in (B_{i,j} \otimes M_{q})^{1}_{+}$ pairwise orthogonal and purely positive, and $\nu_{1},\ldots,\nu_{L} \in [0,1]$ such that for any $\tau \in T(B)$
\begin{eqnarray*}
\textstyle
|\tau(b) - \sum_{l=1}^{L} \nu_{l} \cdot d_{\tau \otimes \mathrm{tr}_{q}}(b_{l}) | & < &  \frac{\epsilon_{i}}{3\cdot M_{i}}.
\end{eqnarray*}
But then, for any $\tau \in T(A)$ we also have 
\begin{eqnarray*}
\textstyle
|(\tau \otimes \tau_{\mathcal{Q}})(\beta_{i}(b)) - \sum_{l=1}^{L} \nu_{l} \cdot d_{\tau \otimes \tau_{\mathcal{Q}} \otimes \mathrm{tr}_{q}}(\beta_{i} \otimes \mathrm{id}_{M_{q}}(b_{l})) | & < & \frac{\epsilon_{i}}{3\cdot M_{i}}.
\end{eqnarray*}
We compute for $\tau \in T(A)$ and $l \in \{1,\ldots,L\}$
\begin{eqnarray*}
d_{\tau \otimes \tau_{\mathcal{Q}} \otimes \mathrm{tr}_{q}}(\beta_{i} \otimes \mathrm{id}_{M_{q}}(b_{l})) & = & d_{\tau \otimes \tau_{\mathcal{Q}} \otimes \mathrm{tr}_{q}}(\kappa_{i}([b_{l}])) \\
& = & \textstyle d_{\tau \otimes \tau_{\mathcal{Q}} \otimes \mathrm{tr}_{q}}(\delta^{-1}(  \frac{\mu_{i,j}}{\lambda_{i,j}} \cdot \gamma ([b_{l} \otimes 1_{\mathcal{Q}}]))) \\
& = & \textstyle   \frac{\mu_{i,j}}{\lambda_{i,j}} \cdot \gamma ([b_{l} \otimes 1_{\mathcal{Q}}]) (T(\iota \otimes \mathrm{id}_{M_{q}})^{-1}(\tau \otimes \mathrm{tr}_{q}) \otimes \tau_{\mathcal{Q}}) \\
& = & \textstyle   \frac{\mu_{i,j}}{\lambda_{i,j}} \cdot d_{T(\iota \otimes \mathrm{id}_{M_{q}})^{-1}(\tau \otimes \mathrm{tr}_{q})}(b_{l})\\
& = & \textstyle   \frac{\mu_{i,j}}{\lambda_{i,j}} \cdot d_{T(\iota)^{-1}(\tau) \otimes \mathrm{tr}_{q}}(b_{l}).
\end{eqnarray*}
It follows that for $\tau \in \partial_{e}T(A)$
\begin{eqnarray*}
\lefteqn{ |(T(\iota)^{-1}(\tau))(b) - (\tau \otimes \tau_{\mathcal{Q}})(\beta_{i}(b)) | }\\
& < & \textstyle | \sum_{l=1}^{L} \nu_{l} \cdot (d_{T(\iota)^{-1}(\tau) \otimes \mathrm{tr}_{q}}(b_{l}) - d_{\tau \otimes \tau_{\mathcal{Q}} \otimes \mathrm{tr}_{q}} (\beta_{i} \otimes \mathrm{id}_{M_{q}}(b_{l}))) | + \frac{2 \epsilon_{i}}{3 M_{i}} \\
& = & \textstyle | \sum_{l=1}^{L} \nu_{l} \cdot (1-\frac{\mu_{i,j}}{\lambda_{i,j}}) \cdot   d_{T(\iota)^{-1}(\tau) \otimes \mathrm{tr}_{q}}(b_{l}) | + \frac{2 \epsilon_{i}}{3 M_{i}} \\
& < & \textstyle |1- \frac{\mu_{i,j}}{\lambda_{i,j}}| \cdot 1+ \frac{2 \epsilon_{i}}{3 M_{i}} \\
& < & \frac{\epsilon_{i}}{M_{i}}.
\end{eqnarray*}
As a consequence, we have 
\[
| (T(\iota)^{-1}(\tau))(b) - (\tau \otimes \tau_{\mathcal{Q}}) (\beta_{i}(b)) | < \epsilon_{i}
\]
for all $b \in (B_{i})^{1}_{+}$ and $\tau \in T(A)$. 

Set
\[
\sigma_{i}:= \psi_{i} \circ \iota: A \to B_{i}
\]
and
\[
\rho_{i}:= B_{i} \stackrel{\beta_{i}}{\longrightarrow} A \otimes \mathcal{Q} \cong A
\]
for $i \in \mathbb{N}$. We estimate
\begin{eqnarray*}
\lefteqn{ \sup_{\tau \in T(A)} \|\tau(\rho_{i} \sigma_{i} (a) - a) \|   } \\
& = & \sup_{\tau \in T(A)} | (\tau \otimes \tau_{\mathcal{Q}}) (\beta_{i} \sigma_{i} (a)) - \tau(a) | \\ 
& \le & \sup_{\tau \in T(A)} | (T(\iota)^{-1}(\tau)) (\sigma_{i} (a)) - \tau(a) | + \epsilon_{i} \\ 
& = & \sup_{\tau \in T(A)} | (T(\iota)^{-1}(\tau)) (\psi_{i} (\iota (a)) - \iota(a)) | + \epsilon_{i} \\ 
& \le & \sup_{\tau \in T(A)} | (T(\iota)^{-1}(\tau)) (\psi_{i} (\iota (a)) - 1_{B_{i}}\iota(a) 1_{B_{i}}) |  \\ 
&& + \sup_{\tau \in T(A)} | (T(\iota)^{-1}(\tau)) (1_{B} -  1_{B_{i}}) |  \\
&& + 2 \|[\iota(a),1_{B_{i}} ] \|\\
&& + \epsilon_{i} \\
& \le & \|\psi_{i}(\iota(a)) - 1_{B_{i}} \iota(a) 1_{B_{i}} \| \\
&& + \sup_{\tau \in T(A)} | (T(\iota)^{-1}(\tau)) (1_{B} -  1_{B_{i}}) |  \\
&& + 2 \|[\iota(a),1_{B_{i}} ] \|\\
&& + \epsilon_{i} \\
& \stackrel{i \to \infty}{\longrightarrow} & 0.
\end{eqnarray*}

We have now verified the hypotheses of Theorem~\ref{A}, whence $A \cong A \otimes \mathcal{Q}$ is TAI.
\end{nproof}
\en

\bn
\label{O}
The following is shown by Strung in \cite{Str:odd-spheres}; this in particular confirms the hypotheses of Theorem~\ref{F} -- see Corollary~\ref{G} below. The Corollary as it stands only covers the examples of \cite{Windsor:finite-ergodic}, but after this work (and that of Strung) was completed, Lin generalized it to cover arbitrary minimal homeomorphisms of odd dimensional spheres; see \cite{Lin:odd-spheres}. Lin's method still uses our classification by embedding technique and also ideas from \cite{Str:odd-spheres}. 

\begin{ntheorem}
Let $\alpha: S^{2n+1} \to S^{2n+1}$, $n \ge 1$, be a minimal homeomorphism which can uniformly be approximated by periodic homeomorphisms as in \cite{Windsor:finite-ergodic}. 

Then, there is a uniquely ergodic minimal homeomorphism $\gamma:X \to X$ on the Cantor set $X$ (in fact, $\gamma$ can be chosen to be an odometer action), such that the product $\alpha \times \gamma: S^{2n+1} \times X \to S^{2n+1} \times X$ is minimal and such that $\mathcal{C}(X \times S^{2n+1}) \rtimes \mathbb{Z}$ is TAI with a uniquely determined state on $\mathrm{K}_{0}$. 
\end{ntheorem}
\en

\bn
\label{G}
\begin{ncor}
Let $\alpha: S^{2n+1} \to S^{2n+1}$, $n \ge 1$, be a minimal homeomorphism, as constructed in \cite{Windsor:finite-ergodic}.

Then, $(\mathcal{C}(S^{2n+1}) \rtimes \mathbb{Z}) \otimes \mathcal{Q}$ is TAI and $\mathcal{C}(S^{2n+1}) \rtimes \mathbb{Z}$ is classifiable in the sense of \cite{Lin:asu-class}. In particular, crossed products of this form are determined by their trace spaces.
\end{ncor}

\begin{nproof}
By Theorem~\ref{O}, there is a uniquely ergodic Cantor minimal action $\gamma:X \to X$ such that $\mathcal{C}(S^{2n+1} \times X) \rtimes \mathbb{Z}$ is TAI with only one state on $\mathrm{K}_{0}$. Let 
\[
\iota: \mathcal{C}(S^{2n+1}) \rtimes \mathbb{Z} \to \mathcal{C}(S^{2n+1} \times X) \rtimes \mathbb{Z}
\]
be the canonical embedding. It is straightforward to check that $T(\iota)$ is a homeomorphism. 

By Theorem~\ref{F}, $(\mathcal{C}(S^{2n+1}) \rtimes \mathbb{Z}) \otimes \mathcal{Q}$ is TAI, hence classifiable (the algebra is in the bootstrap class and satisfies the UCT, cf.\ \cite{Bla:k-theory}). All such crossed products have the same ordered $\mathrm{K}$-theory and only one state on $\mathrm{K}_{0}$, so their Elliott invariants are just distinguished by their trace spaces.
\end{nproof}
\en

\section{Generic classifiability}
\label{section5}

\noindent
We now follow up on a theme from \cite{HirWinZac:Rokhlin-dimension} to show that, in fair generality, classifiability is generically preserved under crossed products by $\mathbb{Z}^{d}$.
We need a notion of Rokhlin dimension for $\mathbb{Z}^{d}$-actions; this is the obvious generalization of \cite[Definition~2.3({c})]{HirWinZac:Rokhlin-dimension} where $\{0,\ldots,p\}$ is replaced by $\{0,\ldots,p\}^{d}$, see \cite{Win:dynamic-dimension} and \cite{Sza:dimnuc} for details.

\bn
\label{L} 
\begin{nprop}
Let $A$ be a separable, simple, unital $\mathrm{C}^{*}$-algebra with $\dim_{\mathrm{nuc}}A< \infty$. Let $\alpha: \mathbb{Z}^{d} \to \mathrm{Aut}(A)$ be an action with finite Rokhlin dimension.

Then, $A \rtimes_{\alpha} \mathbb{Z}^{d}$ is again simple with finite nuclear dimension, and the canonical map $T(A \rtimes_{\alpha} \mathbb{Z}^{d}) \to T(A)^{\alpha}$ is an isomorphism.
\end{nprop}

\begin{nproof}
This is just a variation of \cite[Theorem~4.1]{HirWinZac:Rokhlin-dimension}, see also \cite{Sza:dimnuc}. We omit the details. 
\end{nproof}
\en

\bn
\label{I}
\begin{nprop}
Let $A$ be a separable, unital, $\mathcal{Z}$-stable $\mathrm{C}^{*}$-algebra.

Then, there is a dense $G_{\delta}$ of $\mathbb{Z}^{d}$-actions of $A$ with finite Rokhlin dimension.
\end{nprop}

\begin{nproof}
Apply the technique of \cite[Theorem~3.4]{HirWinZac:Rokhlin-dimension} with $\mathcal{Z}^{\otimes d}$ in place of $\mathcal{Z}$ to end up with a generic set of $\mathbb{Z}^{d}$-actions with Rokhlin dimension at most $2^{d}-1$. (Alternatively, one could also apply \cite[Theorem~3.4]{HirWinZac:Rokhlin-dimension} inductively $d$ times.) 
\end{nproof}
\en

\bn
\label{K}
\begin{ntheorem}
Let $A$ be a $\mathrm{C}^{*}$-algebra with finitely many extremal tracial states and satisfying the conditions of Theorem~\ref{H} (i.e., $A$ is separable, simple, unital, nonelementary, has finite decomposition rank, satisfies the UCT and $\mathrm{K}_{0}(A)$ separates the tracial states of $A$).   

Then, for a dense $G_{\delta}$ of  $\mathbb{Z}^{d}$-actions on $A$, $A \rtimes \mathbb{Z}^{d}$ again has finitely many extremal traces and satisfies the conditions of Theorem~\ref{H}. In particular, for generic sets of $\mathbb{Z}^{d}$-actions, being classifiable by ordered $\mathrm{K}$-theory in the sense of Theorem~\ref{H} and having finitely many extremal tracial states passes to crossed products. 
\end{ntheorem}

\begin{nproof}
By Proposition~\ref{I}, there is a dense $G_{\delta}$ of $\mathbb{Z}^{d}$-actions of $A$ with finite Rokhlin dimension; by Proposition~\ref{L}, every such action gives rise to a crossed product $A \rtimes \mathbb{Z}^{d}$ which is simple, has finite nuclear dimension and with tracial state space being a subspace of $T(A)$. But then $T(A \rtimes \mathbb{Z}^{d})$ is finite dimensional, compact and convex, hence has only finitely many extreme points. These are again separated by projections, whence $A \rtimes \mathbb{Z}^{d}$ has real rank zero. Since $A$ is AH (see Theorem~\ref{H}), by \cite[Theorem~9.3 and its proof]{Lin:crossed-product-AF-embedding}, for each $\tau \in T(A \rtimes \mathbb{Z}^{d})$ there is a unital embedding $\alpha: A \rtimes \mathbb{Z}^{d} \to D$ into a simple, unital, AF algebra with unique tracial state $\delta$ and such that $\tau = \delta \circ \alpha$. Now $(A \rtimes \mathbb{Z}^{d}) \otimes \mathcal{Q}$ is TAF by Theorem~\ref{C}.
\end{nproof}
\en

\bibliographystyle{amsplain}

\begin{thebibliography}{10}

\bibitem{AraPerToms:survey}
P.~Ara, F.~Perera, and A.~S. Toms, \emph{{$\mathrm{K}$-theory for operator
  algebras. Classification of $\mathrm{C}^*$-algebras}}, Aspects of operator
  algebras and applications, Contemp. Math., vol. 534, Amer. Math. Soc,
  Providence RI., 2011, pp.~1--71.

\bibitem{Bla:k-theory}
B.~Blackadar, \emph{{$\mathrm{K}$-Theory for Operator Algebras}}, MSRI
  Monographs, vol.~5, Springer Verlag, Berlin and New York, 1986.

\bibitem{BlaKir:limits}
B.~Blackadar and E.~Kirchberg, \emph{{Generalized inductive limits of
  finite-dimensional $\mathrm{C}^{*}$-algebras}}, Math. Ann. \textbf{307}
  (1997), 343--380.

\bibitem{Con:Thom}
A.~Connes, \emph{{An analogue of the Thom isomorphism for crossed products of a
  $\mathrm{C}^*$-algebra by an action of $\mathbb{R}$}}, Adv. Math. \textbf{39}
  (1981), 31--55.

\bibitem{EllToms:BullAMS}
G.~A. Elliott and A.~S. Toms, \emph{{Regularity properties in the
  classification program for separable amenable {$\mathrm{C}^*$}-algebras}},
  Bull. Amer. Math. Soc. (N.S.) \textbf{45} (2008), no.~2, 229--245.

\bibitem{FatHer:diffmin}
A.~Fathi and M.~R. Herman, \emph{{Existence de diff\'eomorphismes minimaux}},
  Dynamical systems, {V}ol. {I}---{W}arsaw, Soc. Math. France, Paris, 1977,
  pp.~37--59. Ast\'erisque, No. 49.

\bibitem{GMPS:orbit-equivalence-Z2}
T.~Giordano, H.~Matui, I.~F. Putnam, and Ch.~F. Skau, \emph{Orbit equivalence
  for {C}antor minimal {$\Bbb Z^2$}-systems}, J. Amer. Math. Soc. \textbf{21}
  (2008), no.~3, 863--892.

\bibitem{GMPS:orbit-equivalence-Zd}
\bysame, \emph{Orbit equivalence for {C}antor minimal {$\Bbb Z^d$}-systems},
  Invent. Math. \textbf{179} (2010), no.~1, 119--158.

\bibitem{GPS:orbit}
T.~Giordiano, I.~F. Putnam, and C.~F. Skau, \emph{{Topological orbit
  equivalence and $\mathrm{C}^*$-crossed products}}, J. Reine Angew. Math.
  \textbf{469} (1995), 51--111.

\bibitem{HirWinZac:Rokhlin-dimension}
I.~Hirshberg, W.~Winter, and J.~Zacharias, \emph{{Rokhlin dimension and
  $\mathrm{C}^*$-dynamics}}, Comm. Math. Phys. \textbf{335} (2015), 637--670.

\bibitem{JiaSu:Z}
X.~Jiang and H.~Su, \emph{On a simple unital projectionless
  {$\mathrm{C}^*$}-algebra}, Amer. J. Math. \textbf{121} (1999), no.~2,
  359--413.

\bibitem{KirRor:pi2}
E.~Kirchberg and M.~R{\o}rdam, \emph{{Infinite non-simple
  $\mathrm{C}^*$-algebras: absorbing the Cuntz algebra $\mathcal{O}_\infty$}},
  Adv. Math. \textbf{167} (2002), no.~2, 195--264.

\bibitem{KirRor:pi3}
\bysame, \emph{{Purely infinite $\mathrm{C}^*$-algebras: Ideal preserving zero
  homotopies}}, Geom. Funct. Anal. \textbf{15} (2005), no.~2, 377--415.

\bibitem{KirWinter:dr}
E.~Kirchberg and W.~Winter, \emph{{Covering dimension and quasidiagonality}},
  Internat. J. Math. \textbf{15} (2004), 63--85.

\bibitem{Lin:TAFduke}
H.~Lin, \emph{{Classification of simple $\mathrm{C}^*$-algebras of tracial
  topological rank zero}}, Duke Math. J. \textbf{125} (2004), 91--119.

\bibitem{Lin:crossed-product-AF-embedding}
\bysame, \emph{{AF}-embeddings of the crossed products of {AH}-algebras by
  finitely generated abelian groups}, Int. Math. Res. Pap. IMRP (2008), no.~3,
  Art. ID rpn007, 67.

\bibitem{Lin:asu-class}
\bysame, \emph{{Asymptotic unitary equivalence and classification of simple
  amenable $\mathrm{C}^*$-algebras}}, Invent. Math. \textbf{183} (2011), no.~2,
  385--450.

\bibitem{Lin:odd-spheres}
\bysame, \emph{{Minimal dynamical systems on connected odd dimensional
  spaces}}, arXiv preprint math.OA/1404.7034, 2014.

\bibitem{LinNiu:KKlifting}
H.~Lin and Z.~Niu, \emph{Lifting {$\mathrm{KK}$}-elements, asymptotic unitary
  equivalence and classification of simple {$\mathrm{C}^*$}-algebras}, Adv.
  Math. \textbf{219} (2008), no.~5, 1729--1769.

\bibitem{MatuiSato:Comp}
H.~Matui and Y.~Sato, \emph{Strict comparison and {$\mathcal{Z}$}-absorption of
  nuclear {$\mathrm{C}^*$}-algebras}, Acta Math. \textbf{209} (2012), no.~1,
  179--196.

\bibitem{MatSat:dr-UHF}
\bysame, \emph{{Decomposition rank of UHF-absorbing $\mathrm{C}^*$-algebras}},
  Duke Math. J. \textbf{163} (2014), 2687--2708.

\bibitem{Rob:NCCW}
L.~Robert, \emph{Classification of inductive limits of 1-dimensional {NCCW}
  complexes}, Adv. Math. \textbf{231} (2012), no.~5, 2802--2836.

\bibitem{Ror:Z-absorbing}
M.~R{\o}rdam, \emph{The stable and the real rank of {$\mathcal{Z}$}-absorbing
  {$\mathrm{C}^*$}-algebras}, Internat. J. Math. \textbf{15} (2004), no.~10,
  1065--1084.

\bibitem{RorWin:Z-revisited}
M.~R{\o}rdam and W.~Winter, \emph{The {Jiang--Su} algebra revisited}, J. Reine
  Angew. Math. \textbf{642} (2010), 129--155.

\bibitem{Str:thesis}
K.~R. Strung, \emph{{On classification, UHF-stability, and tracial
  approximation of simple nuclear $\mathrm{C}^*$-algebras}}, Ph.D. thesis,
  M\"unster, 2013.

\bibitem{Str:odd-spheres}
\bysame, \emph{{$\mathrm{C}^*$-algebras of minimal dynamical systems of the
  product of a Cantor set and an odd dimensional sphere}}, to appear in J.
  Funct. Anal.; arXiv preprint math.OA/1403.3136, 2014.

\bibitem{StrWin:UHFslicing}
K.~R. Strung and W.~Winter, \emph{{UHF-slicing and classifcation of nuclear
  {$\mathrm{C}^*$}-algebras}}, J. Top. Anal. \textbf{6} (2014), 465--540.

\bibitem{Sza:dimnuc}
G.~Szab\'o, \emph{{The Rokhlin dimension of topological
  $\mathbb{Z}^m$-actions}}, to appear in Proc. London Math. Soc. (3); arXiv
  preprint math.OA/1308.5418, 2013.

\bibitem{TikWin:Z-dr}
A.~Tikuisis and W.~Winter, \emph{{Decomposition rank of $\mathcal{Z}$-stable
  $\mathrm{C}^{*}$-algebras}}, Analysis \& PDE \textbf{7} (2014), 673--700.

\bibitem{TomsWinter:PNAS}
A.~S. Toms and W.~Winter, \emph{{Minimal dynamics and the classification of
  {$\mathrm{C}^*$}-algebras}}, Proc. Natl. Acad. Sci. USA \textbf{106} (2009),
  no.~40, 16942--16943.

\bibitem{TomsWinter:VI}
\bysame, \emph{{The {E}lliott conjecture for {V}illadsen algebras of the first
  type}}, J. Funct. Anal. \textbf{256} (2009), no.~5, 1311--1340.

\bibitem{TomsWinter:minhom}
\bysame, \emph{{Minimal dynamics and $\mathrm{K}$-theoretic rigidity: Elliott's
  conjecture}}, Geom. Funct. Anal. \textbf{23} (2013), no.~1, 467--481.

\bibitem{Windsor:finite-ergodic}
A.~Windsor, \emph{Minimal but not uniquely ergodic diffeomorphisms}, Smooth
  ergodic theory and its applications ({S}eattle, {WA}, 1999), Proc. Sympos.
  Pure Math., vol.~69, Amer. Math. Soc., Providence, RI, 2001, pp.~809--824.

\bibitem{Win:Z-class}
W.~Winter, \emph{On the classification of simple {$\mathcal{Z}$}-stable
  {$\mathrm{C}^{*}$}-algebras with real rank zero and finite decomposition
  rank}, J. London Math. Soc. (2) \textbf{74} (2006), no.~1, 167--183.

\bibitem{Winter:dr-Z-stable}
\bysame, \emph{{Decomposition rank and $\mathcal{Z}$-stability}}, Invent. Math.
  \textbf{179} (2010), no.~2, 229--301.

\bibitem{Winter:dimnuc-Z-stable}
\bysame, \emph{{Nuclear dimension and $\mathcal{Z}$-stability of pure
  $\mathrm{C}^*$-algebras}}, Invent. Math. \textbf{187} (2012), 259--342.

\bibitem{Win:dynamic-dimension}
\bysame, \emph{Dynamic dimension}, in preparation, 2013.

\bibitem{Win:localizingEC}
\bysame, \emph{{Localizing the Elliott conjecture at strongly self-absorbing
  $\mathrm{C}^{*}$-algebras. With an appendix by H.\ Lin.}}, J. Reine Angew.
  Math. \textbf{692} (2014), 193--231.

\bibitem{WinZac:order-zero}
W.~Winter and J.~Zacharias, \emph{Completely positive maps of order zero},
  M\"unster J. Math. \textbf{2} (2009), 311--324.

\bibitem{WinterZac:dimnuc}
\bysame, \emph{{The nuclear dimension of $\mathrm{C}^{*}$-algebras}}, Adv.
  Math. \textbf{224} (2010), 461--498.

\end{thebibliography}

\providecommand{\bysame}{\leavevmode\hbox to3em{\hrulefill}\thinspace}
\providecommand{\MR}{\relax\ifhmode\unskip\space\fi MR }
\providecommand{\MRhref}[2]{%
  \href{http://www.ams.org/mathscinet-getitem?mr=#1}{#2}
}
\providecommand{\href}[2]{#2}

\end{document}